\documentclass[11pt]{amsart}
\usepackage[margin=1.1in]{geometry}
\usepackage{amssymb,booktabs,xcolor,hyperref}
\hypersetup{colorlinks=true,linkcolor=blue!50!black,citecolor=blue!50!black,urlcolor=blue!50!black}

\newtheorem{theorem}{Theorem}
\newtheorem{proposition}[theorem]{Proposition}
\newtheorem{lemma}[theorem]{Lemma}
\theoremstyle{remark}

\newcommand{\R}{\mathbb{R}}
\newcommand{\F}{\mathbb{F}}
\newcommand{\Z}{\mathbb{Z}}

\title[Kissing numbers in dimensions 32 to 39]{New lower bounds for kissing numbers\\ in dimensions 32--35 and 37--39}
\author{Qihong Lin}
\email{linqihong001@qq.com}
\date{October 22, 2026}

\begin{document}

\begin{abstract}
We show that the kissing numbers in dimensions $32$, $33$, $34$, $35$, $37$, $38$ and $39$ are at least $352320$, $370240$,
$390980$, $415308$, $499816$, $605322$ and $768916$, improving the previous bounds $347584$, $364736$, $384836$, $409676$,
$498024$, $591900$ and $763668$. Except in dimension $38$, the configurations are those of Edel, Rains and Sloane with new
binary constant-weight codes of weight $8$ and minimum distance $8$ in the middle layer; each additional word adds $128$
points. In dimension $38$ we add $6855$ antipodal pairs of vectors of the second shell of the Leech lattice to the empty
equator of Kravatskiy's $591612$-point configuration, where $144$ pairs were known; this is an independent set problem in a
regular graph of degree $11178$ on $8386560$ vertices. The codes and lines were found by local search, and all
configurations were checked pair by pair in exact integer arithmetic.
\end{abstract}

\maketitle

\section{Introduction}

The kissing number $\tau_n$ is the largest number of non-overlapping unit balls in $\R^n$ that can touch a common unit
ball. Equivalently, it is the largest size of a \emph{kissing configuration}: a set of nonzero vectors in $\R^n$ any two of
which make an angle of at least $60^\circ$, that is, $\langle x,y\rangle\le\frac12|x||y|$ for $x\ne y$. The exact value is
known only for $n\in\{1,2,3,4,8,24\}$; see~\cite{CS} for background and Cohn's table~\cite{CohnTable} for the current
bounds.

In 1998 Edel, Rains and Sloane~\cite{ERS} gave a simple construction from binary codes that set new records in dimensions
$32$ to $128$. In the form used here it has a dense layer of $\pm1$ vectors from a binary code $D$, a middle layer of $\pm2$
vectors on the words of a constant-weight code, and a layer of scaled roots, and it gives
\begin{equation}\label{eq:ers}
\tau_n\ \ge\ |D|+128\,A(n,8,8)+2n(n-1),
\end{equation}
where $A(n,d,w)$ denotes the largest size of a binary code of length $n$, constant weight $w$ and minimum distance
$d$~\cite{BSSS,BrouwerTable}. With the codes then available, Edel, Rains and Sloane obtained $\tau_{32}\ge276032$ from
$A(32,8,8)\ge1117$ and $\tau_{36}\ge438872$ from $A(36,8,8)\ge2385$, and they expected their bounds to be superseded soon
because the codes were weak. Remarkably, \eqref{eq:ers} and its variants~\cite{YZW} still give the records in most of
dimensions $32$ to $39$, now through much better codes. In 2026 a series of computer searches raised the lower bounds on
$A(n,8,8)$: Rosin~\cite{Rosin}, Echols~\cite{Echols}, Yang et al.~\cite{YZW}, and several contributions recorded in Brouwer's
table~\cite{BrouwerTable}. Dimension $38$ is different: there Kravatskiy~\cite{Kravatskiy} built a $591612$-point
configuration from the minimal vectors of the Leech lattice and Ganzhinov's $1932$-point kissing configuration in
$\R^{14}$~\cite{Ganzhinov}, and Yang et al.\ added $144$ antipodal pairs on its equator, obtaining $591900$~\cite{YZW}.

Our main result is the following.

\begin{theorem}\label{thm:kissing}
$\tau_{32}\ge352320$, $\tau_{33}\ge370240$, $\tau_{34}\ge390980$, $\tau_{35}\ge415308$, $\tau_{37}\ge499816$,
$\tau_{38}\ge605322$ and $\tau_{39}\ge768916$.
\end{theorem}

These improve the previous bounds by $4736$, $5504$, $6144$, $5632$, $1792$, $13422$ and $5248$, respectively; see
Table~\ref{tab:results}. Except for $n=38$, Theorem~\ref{thm:kissing} follows from~\eqref{eq:ers} and new constant-weight
codes.

\begin{theorem}\label{thm:codes}
$A(32,8,8)\ge1713$, $A(33,8,8)\ge1852$, $A(34,8,8)\ge2013$, $A(35,8,8)\ge2202$, $A(37,8,8)\ge2860$,
$A(38,8,8)\ge3088$ and $A(39,8,8)\ge3424$.
\end{theorem}

\begin{table}[h]
\centering
\footnotesize
\setlength{\tabcolsep}{3.5pt}
\begin{tabular}{rrlrrcr}
\toprule
$n$ & previous & source of previous bound & \textbf{new} & gain & here & upper bound \\
\midrule
32 & 347584 & Yang et al.~\cite{YZW}$^{*}$ & \textbf{352320} & 4736 & \S\ref{sec:ers} & 3162316 \\
33 & 364736 & \eqref{eq:ers}, $A(33,8,8)\ge1809$ by Pliny the Liberator~\cite{BrouwerTable} & \textbf{370240} & 5504 & \S\ref{sec:ers} & 4494570 \\
34 & 384836 & \eqref{eq:ers}, $A(34,8,8)\ge1965$ by Pliny the Liberator~\cite{BrouwerTable} & \textbf{390980} & 6144 & \S\ref{sec:ers} & 6422593 \\
35 & 409676 & Yang et al.~\cite{YZW} & \textbf{415308} & 5632 & \S\ref{sec:ers} & 9162403 \\
36 & 484760 & Yang et al.~\cite{YZW} & -- & -- & & 13017098 \\
37 & 498024 & Yang et al.~\cite{YZW} & \textbf{499816} & 1792 & \S\ref{sec:ers} & 18498316 \\
38 & 591900 & Yang et al.~\cite{YZW}, extending Kravatskiy~\cite{Kravatskiy} & \textbf{605322} & 13422 & \S\ref{sec:d38} & 26496684 \\
39 & 763668 & Yang et al.~\cite{YZW} & \textbf{768916} & 5248 & \S\ref{sec:ers} & 37826766 \\
\bottomrule
\end{tabular}
\caption{Lower bounds for $\tau_n$. The previous bounds are those of Cohn's table~\cite{CohnTable} on October 22, 2026; the
upper bounds are from~\cite{LdL}, as listed there. Dimension $36$ is not improved: our best code there has $2742$ words, the
size already known~\cite{BrouwerTable}. $^{*}$The table now lists $352192$ for $n=32$, an earlier version of the present
result~\cite{Data32}.}
\label{tab:results}
\end{table}

\subsection*{What is new}
The constructions themselves are known: the Edel--Rains--Sloane configurations~\cite{ERS}, the dense codes
of~\cite{ChengSloane,ZL,YZW}, Kravatskiy's configuration~\cite{Kravatskiy} and the equatorial lemma of Yang et
al.~\cite{YZW}. What is new is
\begin{itemize}
\item the seven constant-weight codes of Theorem~\ref{thm:codes}, found by an iterated local search combined with large
ruin-and-recreate moves and with a relay that passes lengthened and shortened codes between neighbouring lengths
(Section~\ref{sec:search});
\item a description of the equatorial problem in dimension $38$ as an independent set problem in an explicit regular graph
of degree $11178$ (Lemma~\ref{lem:graph}), whose low degree explains why far more than $144$ lines can be found, and a set of
$6855$ lines (Section~\ref{sec:d38});
\item public certificates and an exact all-pairs verification of every configuration (Section~\ref{sec:data}).
\end{itemize}
Section~\ref{sec:ers} recalls the three-layer construction and gives the codes, Section~\ref{sec:d38} treats dimension $38$,
Section~\ref{sec:search} describes the searches, Section~\ref{sec:data} the verification and the data, and
Section~\ref{sec:remarks} discusses limits and open questions.

\section{Three-layer configurations}\label{sec:ers}

Index the coordinates of $\R^n$ by $\{0,1,\dots,n-1\}$. For a binary word $c$ of length $m\le n$, let $(-1)^c\in\R^n$ have
entries $(-1)^{c_i}$ for $i<m$ and $0$ for $i\ge m$.

\begin{proposition}[Edel--Rains--Sloane~\cite{ERS}]\label{prop:ers}
Let $32\le m\le n$, let $D\subseteq\F_2^m$ have minimum distance at least $m/4$, and let $\mathcal{B}$ be a family of
$8$-subsets of $\{0,\dots,n-1\}$ any two of which meet in at most $4$ points. The following vectors form a kissing
configuration in $\R^n$ of size $|D|+128|\mathcal{B}|+2n(n-1)$:
\begin{enumerate}
\item the dense layer: $(-1)^c$ for $c\in D$ (squared norm $m$);
\item the middle layer: for each $S\in\mathcal{B}$, the $128$ vectors equal to $\pm2$ on $S$, with an even number of minus
signs, and $0$ outside $S$ (squared norm $32$);
\item the root layer: $\pm4e_i\pm4e_j$ for $i<j$ (squared norm $32$).
\end{enumerate}
\end{proposition}

\begin{proof}
Let $x\ne y$ be two of the vectors. If both are dense, $\langle x,y\rangle=m-2d(c,c')\le m/2$. If $x$ is dense and $y$ is
not, $|\langle x,y\rangle|\le16$ (or $\le8$ for a root), and $16\le\frac12\sqrt{32m}$ because $m\ge32$. Two middle vectors
on the same support have distinct even sign patterns, which differ in at least two places, so
$\langle x,y\rangle\le4(8-2\cdot2)=16$; on different supports $S,T$ we have $|\langle x,y\rangle|\le4|S\cap T|\le16$. A middle
vector and a root have inner product at most $2\cdot2\cdot4=16$, and two distinct roots at most $16$. In all cases
$\langle x,y\rangle\le\frac12|x||y|$.
\end{proof}

This is the case of three layers with supports of sizes $m$, $8$ and $2$ of the general construction in~\cite[\S2]{ERS}.
There all layers are scaled to the same norm; we keep integer coordinates instead, which changes nothing since only angles
matter. Dense codes on fewer than $n$ coordinates already appear in~\cite[\S3]{ERS} for $n=36$. A family $\mathcal{B}$ as in
Proposition~\ref{prop:ers} is exactly a binary code of length $n$, weight $8$ and minimum distance $8$. For the dense code
we use:
\begin{itemize}
\item for $32\le n\le37$, the Cheng--Sloane $[32,17,8]$ code~\cite{ChengSloane,Grassl}, as in~\cite{ERS}, so $m=32$ and
$|D|=2^{17}=131072$;
\item for $n=39$, a code with $327680$ words and minimum distance $10$ on all $39$ coordinates, so $m=39$. The bound
$A(39,10)\ge327680$ goes back to Zinov'ev and Litsyn~\cite{ZL}; we use the explicit form of Yang et
al.~\cite[Prop.~5.5 and App.~B]{YZW}, the union of $10$ cosets of a $[39,15,12]$ code obtained by shortening the extended
binary quadratic-residue code of length $48$.
\end{itemize}
Table~\ref{tab:layers} gives the resulting counts. The code of length $38$ does not give a kissing record: with a
$180224$-word dense code of length $38$ and minimum distance $10$~\cite{ZL,YZW} it yields $578300$, much less than the
configuration of Section~\ref{sec:d38}.

\begin{table}[h]
\centering
\footnotesize
\begin{tabular}{rrrrr}
\toprule
$n$ & dense layer & middle layer & roots & total \\
\midrule
32 & 131072 & $128\cdot1713=219264$ & 1984 & \textbf{352320} \\
33 & 131072 & $128\cdot1852=237056$ & 2112 & \textbf{370240} \\
34 & 131072 & $128\cdot2013=257664$ & 2244 & \textbf{390980} \\
35 & 131072 & $128\cdot2202=281856$ & 2380 & \textbf{415308} \\
37 & 131072 & $128\cdot2860=366080$ & 2664 & \textbf{499816} \\
39 & 327680 & $128\cdot3424=438272$ & 2964 & \textbf{768916} \\
\bottomrule
\end{tabular}
\caption{The three-layer configurations of Theorem~\ref{thm:kissing}.}
\label{tab:layers}
\end{table}

The previous configurations of Yang et al.\ in dimensions $35$ to $37$ also replace part of the middle layer by signed
vectors~\cite{YZW}. Ours do not: each new code is large enough that the plain construction already improves the bound.
Table~\ref{tab:codes} lists the codes together with the Johnson bound
$A(n,8,8)\le\lfloor\frac n8\lfloor\frac{n-1}7\lfloor\frac{n-2}6\lfloor\frac{n-3}5\lfloor\frac{n-4}4\rfloor\rfloor\rfloor\rfloor\rfloor$,
which is still more than twice as large in every case.

\begin{table}[h]
\centering
\footnotesize
\setlength{\tabcolsep}{4pt}
\begin{tabular}{rrrrrl}
\toprule
$n$ & previous~\cite{BrouwerTable} & \textbf{new} & Johnson & point degrees & pairs of words meeting in $0,1,2,3,4$ points \\
\midrule
32 & 1676 & \textbf{1713} & 3540 & 388--446 & 117309, 334378, 627547, 205892, 181202 \\
33 & 1809 & \textbf{1852} & 4092 & 185--471 & 142752, 415875, 714783, 249677, 190939 \\
34 & 1965 & \textbf{2013} & 4585 & 220--502 & 178940, 523325, 821479, 299385, 201949 \\
35 & 2157 & \textbf{2202} & 5140 & 251--537 & 229267, 664241, 955016, 359716, 215061 \\
37 & 2846 & \textbf{2860} & 7492 & 266--641 & 453912, 1279255, 1472770, 627258, 255175 \\
38 & 3077 & \textbf{3088} & 8436 & 342--678 & 562756, 1540993, 1682649, 710496, 269434 \\
39 & 3383 & \textbf{3424} & 9286 & 427--734 & 738851, 1956221, 2024566, 842653, 297885 \\
\bottomrule
\end{tabular}
\caption{The codes of Theorem~\ref{thm:codes} ($w=d=8$), compared with Brouwer's table on October 22, 2026. ``Point degrees''
is the range of the number of words containing a given coordinate; the last column is the distribution of $|S\cap T|$ over
all pairs of words.}
\label{tab:codes}
\end{table}

\section{Dimension 38: lines on the equator}\label{sec:d38}

We use the integer model of the Leech lattice $\Lambda$ with minimum squared norm $32$~\cite{CS}: the vectors
$v\in\Z^{24}$ whose entries all have the same parity $e\in\{0,1\}$, such that $\bigl((v_i-e)/2\bmod2\bigr)_i$ is a word of
the extended Golay code, and $\sum_iv_i\equiv4e\pmod 8$. We take the Golay code with the generator matrix of Kravatskiy's
package~\cite{Kravatskiy}, so that our vectors and his configuration live in the same lattice. Inner products in $\Lambda$
are multiples of $8$. Let $\Lambda_{32}$ and $\Lambda_{48}$ be the shells of squared norm $32$ and $48$; they have $196560$
and $16773120$ vectors, that is, $98280$ and $8386560$ lines.

\subsection*{Kravatskiy's configuration}
Write $\R^{38}=\R^{24}\oplus\R^{14}$. Let $W\subset\R^{14}$ be Ganzhinov's $1932$-point kissing
configuration~\cite{Ganzhinov}, in the coordinates of the data accompanying~\cite{CohnTable}, scaled to squared norm $8$, so that $\langle t,t'\rangle\le4$ for distinct $t,t'\in W$.
Kravatskiy~\cite{Kravatskiy} partitioned $W$ into $644$ triples $\{t_{k,1},t_{k,2},t_{k,3}\}$ with $t_{k,1}+t_{k,2}+t_{k,3}=0$,
and the $98280$ lines of $\Lambda_{32}$ into $644$ classes such that any two lines in the same class satisfy
$|\langle u,u'\rangle|\le8$. His configuration consists of the $6\cdot98280=589680$ \emph{cap points}
$(\pm u,\sqrt2\,t_{k,i})$, for each line $\pm u$ of class $k$ and $i=1,2,3$, and of the $1932$ \emph{axis points}
$(0,\sqrt6\,Rt)$, $t\in W$, where $R$ is a rotation of $\R^{14}$ with rational entries such that $\langle Rt,t'\rangle\le\sqrt{48}$
for all $t,t'\in W$. All points have squared norm $48$. Two cap points have inner product
$\langle u,u'\rangle+2\langle t,t'\rangle\le24$, a cap point and an axis point have inner product
$\sqrt{12}\,\langle Rt',t\rangle\le24$, and two axis points have inner product $6\langle t,t'\rangle\le24$; so this is a
kissing configuration with $591612$ points. Its equator $\R^{24}\oplus0$ contains no point.

\begin{lemma}[Yang et al.~{\cite[Lemma~7.1, Cor.~7.3]{YZW}}]\label{lem:equator}
Let $V\subset\Lambda_{48}$ contain at most one vector from each line, with $|\langle v,w\rangle|\le24$ for all distinct
$v,w\in V$. Then Kravatskiy's configuration together with the $2|V|$ points $(\pm v,0)$, $v\in V$, is a kissing configuration
with $591612+2|V|$ points.
\end{lemma}

\begin{proof}
We repeat the short argument. A new point is orthogonal to every axis point. For a new point $(v,0)$ and a cap point
$(u,\sqrt2\,t)$, both of squared norm $48$, the cosine is $\langle v,u\rangle/48$; since $v-u$ is a nonzero vector of
$\Lambda$, $80-2\langle v,u\rangle=|v-u|^2\ge32$, so $\langle v,u\rangle\le24$. Two new points have cosine
$\pm\langle v,w\rangle/48$, with $|\langle v,w\rangle|\le24$ by assumption.
\end{proof}

So $\tau_{38}\ge591612+2\alpha(\Gamma)$, where $\alpha(\Gamma)$ is the independence number of the graph $\Gamma$ on the
$8386560$ lines of $\Lambda_{48}$, two lines being adjacent when $|\langle v,w\rangle|>24$. Yang et al.\ found an independent
set of $144$ lines~\cite[Thm.~7.5]{YZW}. They also observed that $\Lambda_{48}$ is the only shell for which
compatibility with the caps is automatic, and that no vector of $\Lambda_{64}$ can be added~\cite[Cor.~7.3,
Prop.~7.7]{YZW}.

\begin{lemma}\label{lem:graph}
Two lines $\pm v\ne\pm w$ of $\Lambda_{48}$ are adjacent in $\Gamma$ if and only if $w=\pm(v-u)$ for some $u\in\Lambda_{32}$
with $\langle v,u\rangle=16$. The graph $\Gamma$ is regular of degree $11178$.
\end{lemma}

\begin{proof}
If $v\ne\pm w$ then $96-2\langle v,w\rangle=|v-w|^2\ge32$, so $\langle v,w\rangle\le32$, and likewise
$\langle v,w\rangle\ge-32$. As inner products are multiples of $8$, adjacency means $\langle v,w\rangle=\pm32$. Replacing $w$
by $-w$ if necessary, $\langle v,w\rangle=32$; then $u=v-w$ has squared norm $32$ and $\langle v,u\rangle=16$. Conversely, if
$u\in\Lambda_{32}$ and $\langle v,u\rangle=16$, then $w=v-u$ has squared norm $48$ and $\langle v,w\rangle=32$.

Distinct $u$ give distinct $w=v-u$, and $v-u=-(v-u')$ would give $\langle v,u+u'\rangle=\langle v,2v\rangle=96\ne32$; so the
degree of $v$ is the number $a_2$ of $u\in\Lambda_{32}$ with $\langle v,u\rangle=16$. Let $a_k$ be the number of
$u\in\Lambda_{32}$ with $\langle v,u\rangle=8k$. Then $a_{-k}=a_k$, and $a_k=0$ for $|k|\ge4$ by the proof of
Lemma~\ref{lem:equator}. Since $\Lambda_{32}$ is a spherical $11$-design~\cite{CS},
\[
\sum_k a_k\,(8k)^{2j}\ =\ 196560\cdot(32\cdot48)^j\,\frac{1\cdot3\cdots(2j-1)}{24\cdot26\cdots(22+2j)}\qquad(j=0,1,2,3).
\]
These four equations give $(a_0,a_1,a_2,a_3)=(75900,48600,11178,552)$.
\end{proof}

The value $a_3=552$ is the number of minimal vectors tangent to an added line, also computed in~\cite[Prop.~7.4]{YZW}. By
Lemma~\ref{lem:graph}, a greedy independent set has at least $\lceil8386560/11179\rceil=751$ lines, and in practice randomized
greedy finds about $5300$. With the local search of Section~\ref{sec:eqsearch} we found $6855$ lines, so
\[
\tau_{38}\ \ge\ 591612+2\cdot6855\ =\ 605322 .
\]
In the other direction, the $2|V|$ new points by themselves form a kissing configuration in $\R^{24}$, so
$|V|\le\tau_{24}/2=98280$.

\section{The search}\label{sec:search}

\subsection{Formulation}
Let $G_n$ be the graph on all $8$-subsets of $\{0,\dots,n-1\}$, two subsets being adjacent when they share at least $5$
points. Then $A(n,8,8)$ is the independence number of $G_n$. The graph has $\binom n8$ vertices ($10518300$ for $n=32$,
$61523748$ for $n=39$) and is regular of degree $\sum_{k=5}^{7}\binom8k\binom{n-8}{8-k}$ ($121264$ for $n=32$). It is never
stored: vertices are indexed by their rank in the combinatorial number system and neighbours are enumerated when needed.
The searches only propose candidates; the results rest on the final codes and lines, which are checked as described in
Section~\ref{sec:data}, not on any property of the search.

\subsection{Iterated local search}
We follow the iterated local search of Andrade, Resende and Werneck~\cite{ARW}. For the current independent set $I$ we keep,
for every vertex $v\notin I$, the number $\kappa(v)=|N(v)\cap I|$ (one byte per vertex); vertices with $\kappa(v)\le6$ are
kept in buckets, and a vertex with $\kappa(v)=1$ records its unique neighbour in $I$. The \emph{local search} inserts free
vertices ($\kappa=0$) in random order and applies $(1,2)$-swaps (a member $s$ that is the only neighbour in $I$ of two
non-adjacent vertices is replaced by both) until neither move applies. The \emph{perturbation} chooses a level
$b\in\{1,2,3,4\}$ with probabilities $(p_1,\dots,p_4)$, samples four vertices with $\kappa=b$, forces the least recently
touched one into $I$ (deleting its $b$ neighbours) and runs the local search. If the set shrank by $\delta$ and is $\delta^*$
below the best size seen, the result is kept with probability $1/(1+c\,\delta\,\delta^*)$ and otherwise undone. We used the
parameter sets
\[
\begin{aligned}
(p_1,p_2,p_3,p_4;\,c)\in\{&(0.45,0.35,0.15,0.05;\,1),\ (0.65,0.25,0.08,0.02;\,2),\\
&(0.45,0.35,0.15,0.05;\,4),\ (0.30,0.35,0.25,0.10;\,2)\}.
\end{aligned}
\]

\subsection{Ruin and recreate}
Started from the $1676$-word code of~\cite{YZW}, the iterated local search reached $1695$ words for $n=32$ within an hour on
a laptop and $1696$ on four servers, and then stalled. Following the ruin-and-recreate principle~\cite{RR}, we then repeatedly
delete a large structured part of the best code and rebuild it with the local search, keeping the result only if it is strictly
larger. The ruin step deletes one of: all words meeting a random $8$-subset in at least $2$ points (about $67\%$ of the words
for $n=32$) or in at least $3$ points (about $31\%$); all words meeting a chosen word in at least $2$ points; or a random
$30$--$70\%$ of the words. The first such cycle kept $622$ of the $1696$ words and rebuilt them to $1706$ within three minutes;
further cycles gave $1707$, $1710$, $1711$, $1712$ and $1713$.

\subsection{Explicit candidate pools}
In a variant, each cycle selects about $30000$ vertices of $G_n$ with the fewest neighbours in the current code, builds the
induced subgraph on the code and these candidates explicitly (average degree about $400$), runs the iterated local search on
this much smaller graph, and then rebuilds the pool.

\subsection{Relay between lengths}\label{sec:relay}
A code of length $n-k$ becomes a code of length $n$ by adding $k$ unused coordinates, and the words of a code of length
$n+1$ that avoid a given coordinate form a code of length $n$. Echols~\cite{Echols} used such lengthened and shortened codes
as seeds for a tabu search. They also turned out to be much better starting points for our searches than the previous record
codes: for $n=35$, searches started from the previous record code stalled below it, while searches started from our codes of
lengths $33$ and $32$, with two and three new coordinates and a greedy completion, passed it. We therefore automated this as a
relay: every $20$ minutes, each new best code was lengthened and shortened and passed to a third of the searches in the
neighbouring lengths. The first starting points were the best published codes~\cite{YZW,BrouwerTable}. The codes of lengths
$32$ and $38$ in Theorem~\ref{thm:codes} were found before the candidate pools and the relay were introduced; the other codes
come from runs that used them.

\subsection{Equatorial lines}\label{sec:eqsearch}
For $\Gamma$ we used the same iterated local search. By Lemma~\ref{lem:graph} the neighbours of a line $v$ are obtained from
one pass over the $196560$ vectors $u\in\Lambda_{32}$, keeping those with $\langle v,u\rangle=16$, followed by a hash-table
lookup of $v-u$; a query takes about a millisecond. Starting from a randomized greedy set of $5291$ lines, the local search
reached $6623$ lines in three minutes, and a ruin-and-recreate loop that deletes a random $10$--$40\%$ of the lines reached
$6855$.

\subsection{Resources}
All searches ran for about a day and a half on four servers with $760$ cores in total and a laptop.

\section{Verification and data}\label{sec:data}

Every claim is decided by exact integer arithmetic on the published data, and the arguments of Sections~\ref{sec:ers}
and~\ref{sec:d38} reduce each bound to finitely many checks:
\begin{enumerate}
\item for $n\ne38$: the dense code has the stated size and minimum distance (enumeration of all words), and the code of
Theorem~\ref{thm:codes} consists of distinct $8$-subsets meeting pairwise in at most $4$ points;
\item for $n=38$: the $6855$ vectors lie in $\Lambda_{48}$ and satisfy $|\langle v,w\rangle|\le24$ pairwise, and Kravatskiy's
configuration has the stated properties.
\end{enumerate}
Item~(2) is checked by a program written from scratch, which also confirms $|\langle v,u\rangle|\le24$ for all $196560$
vectors $u\in\Lambda_{32}$, and by the verifier of Yang et al.~\cite{YZW}, extended to any number of equatorial lines, which
rebuilds Kravatskiy's configuration from its data and checks its classes, triples and rotation.

Independently of all these arguments, a brute-force program, separate from the search code, builds every point of every
configuration as an integer vector and tests every pair $x\ne y$ for
\[
\langle x,y\rangle\le0\quad\text{or}\quad 4\langle x,y\rangle^2\le|x|^2|y|^2,
\]
which is exactly $\langle x,y\rangle\le\frac12|x||y|$ and also excludes repeated or parallel points. In dimension $38$ the
points are written as integer vectors in $\Z^{52}$ through the isometry $(x,y)\mapsto(x,y/\sqrt2,y/\sqrt2)$ of
$\R^{24}\oplus\R^{14}$, with each axis point scaled by an integer; its entries are of order $10^{17}$, and products involving
axis points are computed with $128$-bit integers. The program checks $6.21\cdot10^{10}$, $6.85\cdot10^{10}$,
$7.64\cdot10^{10}$, $8.62\cdot10^{10}$, $1.25\cdot10^{11}$, $1.83\cdot10^{11}$ and $2.96\cdot10^{11}$ pairs for
$n=32,33,34,35,37,38,39$, finds no violation, and finds that the largest cosine is exactly $\frac12$ in each case. All checks
together take about five minutes on a $15$-core laptop.

The certificates and programs are archived on Zenodo: dimension $32$ in~\cite{Data32} and dimensions $33$--$39$
in~\cite{Data33}, together with the base data of Kravatskiy's configuration and the dense code of length $39$ with their
licences. Codes are stored as one $0/1$ string per word (leftmost character $=$ coordinate $0$), the equatorial lines as
$24$ integers per line in the coordinates of the Golay basis above, and each archive has a script that runs all checks.
Table~\ref{tab:sha} gives the SHA-256 hashes of the main files.

\begin{table}[h]
\centering
\scriptsize
\begin{tabular}{ll}
\toprule
file & SHA-256 \\
\midrule
\texttt{supports\_1713.txt} ($n=32$) & \texttt{539972293fbfc671477cea1805a5bf70a2208a4c1e60ec122347b5ef8a23c807} \\
\texttt{d33/supports\_1852.txt} & \texttt{bb3b91ef3a81aa8028f6d170a4743b5788120b7787a5956aaa2c32e5c1ed51df} \\
\texttt{d34/supports\_2013.txt} & \texttt{7a385086eca6429dfe439a73d5d02cd73fa47a98e5b749d0b595dbd623765559} \\
\texttt{d35/supports\_2202.txt} & \texttt{f1049f812373dfbf9473189b4cb5ef8b39bd13c20a9476d881aa3d61f2e448ca} \\
\texttt{d37/supports\_2860.txt} & \texttt{7eef8762bc47991923b0883d557b8cea05027cf31425e6084921494157430610} \\
\texttt{d38/supports\_3088.txt} & \texttt{ccb530f90e914df84cbc8d19a71d30c79146d23a2e8e86e7a813d4d92a5a6a8b} \\
\texttt{d38/equator\_lines\_6855.txt} & \texttt{bd64fdea42bc973bc428f2a60daaf21da5bab33f0d26b0c261a7d70b732a86c8} \\
\texttt{d39/supports\_3424.txt} & \texttt{5e750c76c21aa9e350784624d367cc724a392eb8cb006826613d8dc82321d1cf} \\
\bottomrule
\end{tabular}
\caption{SHA-256 hashes of the certificates in~\cite{Data32,Data33}.}
\label{tab:sha}
\end{table}

\section{Remarks and open questions}\label{sec:remarks}

\emph{Codes.} Each additional word improves~\eqref{eq:ers} by $128$, and the Johnson bound still leaves a factor of more than
two (Table~\ref{tab:codes}); for instance, it limits~\eqref{eq:ers} to $586176$ for $n=32$. We expect that all the codes of
Theorem~\ref{thm:codes} can be improved. Several ideas did not help for $n=32$: merging two or more good codes (solved
exactly as a maximum independent set in a bipartite conflict graph), replacing complete sign groups of the middle layer by
partial sign groups on other supports, and starting from codes with a prescribed symmetry group ($\mathbb{Z}_2^5$ or
$\mathbb{Z}_{31}$) or from random greedy codes.

\emph{Dimension 36.} To beat $484760$ with~\eqref{eq:ers} one needs $A(36,8,8)\ge2744$. The best known code, with $2742$
words, consists of a highly symmetric $2739$-word code and three further words~\cite{BrouwerTable,YZW}. Hundreds of
independent searches from different starting points stopped at $2742$; an exact optimization over the union of the words of
many such codes also gives $2742$; and deleting large structured parts of the code and refilling them exactly gave no gain.
This does not exclude larger codes, but a different construction may be needed.

\emph{Dimension 38.} The independence number $\alpha(\Gamma)$ is only known to lie between $6855$ and $98280$; each further
line adds two points. A spectral or semidefinite bound for $\alpha(\Gamma)$ would show how much room is left. For comparison,
a three-layer configuration with the $180224$-word dense code would need $A(38,8,8)\ge3300$ to reach $605322$.

\section*{Use of artificial intelligence}
This work relied heavily on an AI coding assistant, Claude (model Claude Opus 5.5, Anthropic), used through the Claude Code
tool, at times with several instances working in parallel. For the searches, the assistant proposed and implemented the
ruin-and-recreate loop, the explicit candidate pools and the local search for equatorial lines, implemented the remaining
searches of Section~\ref{sec:search}, and ran all of them on the author's computers. The relay between code lengths grew out
of the author's suggestion to pass improvements on to neighbouring lengths; the assistant turned it into the automatic relay
of Section~\ref{sec:relay}. The assistant also found and proved Lemma~\ref{lem:graph}, checked the arguments taken
from~\cite{ERS,YZW}, wrote the verification programs, checked the literature and the previous records, prepared the data
archives and drafted this paper. The author chose the problems and the dimensions to attack, directed the work and the
allocation of computing resources, and reviewed the results and the text.

The bounds do not depend on trusting the assistant: they rest on the published certificates, which can be checked with the
programs described in Section~\ref{sec:data} or with independent ones. The author takes full responsibility for the content
of this paper.

\subsection*{Acknowledgments}
I thank A.~Kravatskiy and the authors of~\cite{YZW} for making their configurations, data and verification programs public;
the results in dimensions $38$ and $39$ build directly on them.

\end{document}